\documentclass[11pt]{article}
\usepackage[auth-sc]{authblk}
\usepackage{amssymb}
\usepackage{graphicx}
\usepackage{latexsym}
\usepackage{amsmath,amscd}
\usepackage{bm}
\def\part#1{\frac{\partial\phantom{q}}{\partial#1}}

\newenvironment{rmk}{\begin{trivlist}\item[]{\bf Remark:} }
{\end{trivlist}}
\newenvironment{rmks}{\begin{trivlist}\item[]{\bf Remarks:} }
{\end{trivlist}}

\newenvironment{prf}{\begin{trivlist}\item[]{\bf Proof:} }
{\hfill $\Box$ \end{trivlist}}

\newtheorem{thm}{Theorem}

\newtheorem{prp}[thm]{Proposition}

\newcommand{\lie}[1]{\mathfrak{#1}}

\def\Pic{\mathop{\rm Pic}\nolimits}
\def\Hom{\mathop{\rm Hom}\nolimits}

\def\ker{\mathop{\rm ker}\nolimits}
\def\coker{\mathop{\rm coker}\nolimits}

\def\deg{\mathop{\rm deg}\nolimits}

\def\rk{\mathop{\rm rk}\nolimits}

\def\Nm{\mathop{\rm Nm}\nolimits}

\newcommand{\R}{\mathbf{R}}
\newcommand{\C}{\mathbf{C}}
\newcommand{\K}{\mathbf{H}}
\newcommand{\Z}{\mathbf{Z}}

\newcommand{\PP}{{\mathbf {\rm P}}}

\title{ Higgs bundles and characteristic classes}
 \author{Nigel Hitchin}

 \affil  {Mathematical Institute,
Radcliffe Observatory Quarter,
Woodstock Road,
Oxford, OX2 6GG}

\begin{document}
 \maketitle
 \thispagestyle{empty}

\let\oldthefootnote\thefootnote
\renewcommand{\thefootnote}{\fnsymbol{footnote}}
\footnotetext{hitchin@maths.ox.ac.uk}
\let\thefootnote\oldthefootnote

\centerline{\it Dedicated to the memory of Friedrich Hirzebruch} 
\section{Introduction}
Sixty years ago Hirzebruch observed how the vanishing of the Stiefel-Whitney class $w_2$ led to integrality  of the $\hat A$-genus of an algebraic variety \cite{Hirz1}. This was one motivation for  the Atiyah-Singer index theorem but also for my own thesis about Dirac operators and K\"ahler manifolds. Indeed the interaction between topology and algebraic geometry which he developed has been a constant theme in virtually all my work. 

 This article  is also about $w_2$, characteristic classes  and algebraic geometry, but in a rather different context. We consider a compact oriented surface $\Sigma$ of genus $g$, a noncompact real Lie group $G^r$ and the character variety $\Hom(\pi_1(\Sigma),G^r)/G^r$, equivalently the moduli  space of flat $G^r$-connections on $\Sigma$. If $U$ is the maximal compact subgroup of $G^r$ then each principal  $G^r$-bundle  has  a characteristic class in $H^2(\Sigma,\pi_1(U))$ which helps to determine in which  connected component of the character variety it lies. 

A well-known example is  the case $G^r=SL(2,\R)$ where  we have $U=SO(2)$ and a class $c\in H^2(\Sigma,\pi_1(SO(2)))\cong \Z$.  It satisfies the Milnor-Wood inequality $\vert c\vert \le 2g-2$ and there is one component for each $c$ for which strict inequality holds, but  when  $\vert c\vert=2g-2$ there are $2^{2g}$ connected components, each a copy of Teichm\"uller space. 

Here we shall consider the groups $SL(n,\R)$ and $Sp(2m,\R)$, higher dimensional generalizations of $SL(2,\R)=Sp(2,\R)$. In the first case we have a characteristic class in $H^2(\Sigma,\pi_1(SO(n)))$ which, for $n>2$,  is a Stiefel-Whitney class $w_2\in \Z_2$ and in the second a Chern class $c_1$ in $H^2(\Sigma,\pi_1(U(m)))\cong \Z$.

We approach this question using the moduli space of Higgs bundles. For a complex group $G^c$ and a complex structure on $\Sigma$, gauge-theoretic equations enable us to describe the $G^c$-character variety in terms of a holomorphic principal bundle and a holomorphic section $\Phi$ of the associated bundle $\lie{g}\otimes K$ where $K$ is the canonical bundle. This moduli space has the  structure of a completely integrable Hamiltonian system -- a proper map to an affine space, whose generic fibre is an abelian variety. The $G^r$-character variety is realized as the fixed point set of a holomorphic involution and for the real groups in question the involution acts trivially on the base and its fixed points can be identified with the elements of order $2$ in the  abelian variety. The main problem we address is to evaluate the characteristic class as a function on this $\Z_2$-vector space. For $SL(n,\R)$, $w_2$ is a quadratic function related to the mod 2 index theorem and for $Sp(2m,\R)$ the characteristic class is determined by the orbit of an action of a symmetric group under its permutation representation over $\Z_2$. 

The above results provide the background for testing the predictions of mirror symmetry for the hyperk\"ahler Higgs bundle moduli space, and this we approach in the final section of this paper. A component of the $G^r$-character variety is an example of  what is known as a BAA-brane. The SYZ approach to mirror symmetry says that its mirror should be a BBB-brane.  In this context, a BBB-brane is a hyperholomorphic bundle over a hyperk\"ahler submanifold. Now it is known that mirror symmetry for these moduli spaces is closely related to Langlands duality and  the duality of the abelian varieties in the integrable system. How this works for the two real forms above is still a mystery, but we shall indicate a conjectural mirror for the real form $G^r=U(m,m)\subset GL(2m,\C)$. There is again an integral characteristic class here and our candidate for the mirror is, for each allowable value, a hyperholomorphic bundle over the moduli space of $Sp(2m,\C)$-Higgs bundles considered as  a hyperk\"ahler submanifold of the $GL(2m,\C)$-moduli space. It is a trivial line bundle only for the components of the moduli space where the characteristic class takes its  maximum   absolute value. 
\section{Higgs bundles}
We summarize here the basic facts about Higgs bundles  \cite{Hit1},\cite{Sim}. A crucial feature is the hyperk\"ahler structure which provides the non-holomorphic isomorphism  between the character variety and the moduli space of Higgs  bundles. This arises from an infinite-dimensional  quotient construction.

Let $\Sigma$ be a compact oriented Riemann surface of genus $g>1$ and $P$ a principal bundle for the compact real form  $G$ of a complex semi-simple Lie group $G^c$. The space of connections on $P$ is an affine space ${\mathcal A}$ with group of translations $\Omega^1(\Sigma,\lie{g})$ and a symplectic form given by integrating $B(a\wedge b)$ where $B$ is an invariant metric on $G$. The complex structure on $\Sigma$ gives this the structure of an infinite-dimensional flat K\"ahler manifold with complex tangent space $\Omega^{0,1}(\Sigma,\lie{g})$. The group ${\mathcal G}$ of gauge transformations acts isometrically. The cotangent bundle $T^*{\mathcal A}={\mathcal A}\times \Omega^{1,0}(\Sigma,\lie{g})$ is a flat {\it hyperk\"ahler} manifold. This is  a metric compatible with complex structures $I,J,K$ satisfying the relations of the quaternions.  The induced gauge group action has a hyperk\"ahler moment map, which applied to $(A,\Phi)\in {\mathcal A}\times \Omega^{1,0}(\Sigma,\lie{g})$ is 
\begin{equation}
\mu(A,\Phi)=(F_A+[\Phi,\Phi^*], \bar\partial _A\Phi)
\label{mom}
\end{equation}
where $F_A$ is the curvature and $\Phi\mapsto \Phi^*$ for a general group is the antiholomorphic involution coming from  the compact real form. 

The zero set gives firstly $\bar\partial _A\Phi=0$, so the Higgs field $\Phi$ is a holomorphic section of $\lie{g}\otimes K$, and secondly the equation  $F_A+[\Phi,\Phi^*]=0$ which is  equivalent to a stability condition. The quotient of this zero set by ${\mathcal G}$ is the moduli space of pairs $(A,\Phi)$, and it has an induced  hyperk\"ahler structure.  If $I$ denotes the complex structure of pairs $(A,\Phi)$, then $J,K$ are the complex structures for the $G^c$-connections $\nabla_A+\Phi+\Phi^*, \nabla_A+i\Phi-i\Phi^*$ respectively. Setting  (\ref{mom}) to zero shows that these are  flat connections and by considering the holonomy, the moduli space with complex structure $J$ or $K$ can be identified with the $G^c$-character variety.  

The integrable system for  $G^c=GL(n,\C)$ is defined by the characteristic polynomial of the Higgs field  in its defining representation: $\det(x-\Phi)=x^n+a_1x^{n-1}+\dots+a_n$ where $a_i\in H^0(\Sigma,K^i)$. This maps the moduli space ${\mathcal M}$  to the vector space $\bigoplus_{1\le i\le n}H^0(\Sigma,K^i)$: it is proper  and the functions defined by it Poisson-commute with respect to the natural symplectic structure. The generic fibre is then a complex torus but it can be identified with the Jacobian of the curve with equation  $\det(x-\Phi)=0$. For linear groups in general the fibre is a distinguished abelian subvariety of the Jacobian.

To obtain the character variety for the real form $G^r$ in the Higgs bundle realization, we need the connection $A$ to have holonomy in $U$, the maximal compact subgroup of $G^r$ and, with $\lie{g}=\lie{u}\oplus \lie{m}$, the Higgs field must lie in $H^0(\Sigma,\lie{m}\otimes K)$. So for $G^r=SL(n,\R)$, we have $U=SO(n)$ and the Higgs bundle is defined by a rank $n$ holomorphic vector bundle $V$ with an orthogonal structure and $\Lambda^n V$ trivial. The Higgs field    $\Phi$ must then be  symmetric with respect to this inner product \cite{Hit3}. The characteristic class here is $w_2\in \Z_2$: the obstruction to lifting the $SO(n,\C)$-frame bundle to $Spin(n,\C)$.

For $G^r=Sp(2m,\R)$,  the maximal compact is $U=U(m)$ and this means the  vector bundle $V=W\oplus W^*$ with $W$ a rank $m$ vector bundle \cite{GGM},  the pairing between $W$ and $W^*$ defining the symplectic structure.  Here the Higgs field has the off-diagonal form $\Phi(w,\xi)=(\beta(\xi),\gamma(w))$ where $\beta:W^*\rightarrow W\otimes K$ and $\gamma:W\rightarrow W^*\otimes K$ are symmetric. The characteristic class in this case is $c_1(W)\in \Z$. 

\section{The canonical section}\label{can}
 As shown in \cite{Hit3}, there are  canonical sections of the integrable system, each point of which   gives a Higgs bundle for the split real form $G^r$. In particular, this gives a distinguished point in the generic fibre which we can regard as the identity element in an abelian variety. We spell this out next in our two cases which are indeed split real forms.
 
 For $SL(n,\R)$ and $n=2m+1$ the vector bundle is given by 
 $$V=K^{-m}\oplus K^{1-m}\oplus \dots \oplus K^{m}$$
 and for $n=2m$ by
 $$V=K^{-(2m-1)/2}\oplus K^{1-(2m-1)/2}\oplus \dots \oplus K^{(2m-1)/2}$$
 where for $n$ odd we have to choose a square root of the canonical bundle $K$ (a {\it theta characteristic} in classical terms, or a spin structure \cite{MFA1} in the language of topology). The pairing of $K^{\pm \ell}$ or $K^{\pm \ell/2}$ defines an orthogonal structure on $V$ and $\Lambda^nV$ is trivial so it has structure group $SO(n,\C)$. 
 
 The subbundle $K^{1/2}\oplus \dots \oplus K^{(2m-1)/2}$ when $n=2m$ or $K\oplus \dots \oplus K^m$ when $n=2m+1$  is maximal isotropic and a spin structure for $V$ is defined by a holomorphic square root of the top exterior power of a maximal isotropic subbundle.  This in the two cases is
 $K^{m^2/2}$ and $K^{m(m+1)/2}$. These have Chern classes $m^2(g-1)$ and $m(m+1)(g-1)$. The latter is even and so in odd dimensions $w_2=c_1$ mod $2$ $=0$. When $n=2m$,   $w_2=0$ if $g$ is odd and if $g$ is even $w_2=m$ mod $2$.
 
 The Higgs field must be symmetric with respect to this orthogonal structure. We set:
 
  \begin{equation}
 \Phi=\begin{pmatrix}0 & 1 & 0 &\dots & & 0\\
a_2 & 0 & 1 & \dots & & 0\\
a_3 & a_2 & 0 & 1& \dots & 0\\
\vdots & & &\ddots & &\vdots \\
a_{n-1}& & & & \ddots & 1\\
a_n& a_{n-1} & \dots & a_3 & a_2& 0
\end{pmatrix}
 \label{sym}
 \end{equation}
 where $a_i\in H^0(\Sigma, K^i)$. 
 
 \begin{rmk} In \cite{Hit3} p.456 it was claimed that this is conjugate to the companion matrix of the polynomial $x^n+a_2x^{n-2}+\dots+a_n$ which is incorrect. However, the coefficients of the characteristic polynomial are universal polynomials in the $a_i$ and can be thought of as simply changing the basis of invariant polynomials on $\lie{sl}(n)$ which define the fibration.  The actual characteristic polynomial can be viewed as follows (see \cite{Tr}). Set $p(x)= 1-\lambda x+ a_2x^2+\dots + a_nx^n$, then $p(x)$ and $x^n$ have no common factor so there are unique polynomials $a(x), b(x)$ of degree $\le (n-1)$ such that $a(x)p(x)+b(x)x^n=1$ .  Then $b(0)=\det(\lambda-\Phi)$. 
  \end{rmk} 
  
  For $Sp(2m,\R)$ the vector bundle  is given by 
   $$V=K^{-(2m-1)/2}\oplus K^{1-(2m-1)/2}\oplus \dots \oplus K^{(2m-1)/2}$$
  where now we use the pairing between $K^{\pm \ell/2}$ to define a symplectic structure. Putting $W= K^{(2m-1)/2}\oplus K^{(2m-1)/2-2}\oplus \dots \oplus K^{-(2m-3)/2}$ gives the form $V=W\oplus W^*$ above. Then
  $c_1(W)=m(g-1)$.

   We need for the Higgs field sections $\beta,\gamma$ of $S^2W\otimes K, S^2W^*\otimes K$ respectively. Since $W=W^*\otimes K$, we set $\gamma =1$.  Then we take for $\Phi$ the  matrix of the form 
   $$\Phi=\begin{pmatrix} 0 & 1\\
                                 A & 0
                                 \end{pmatrix}$$
  where 
     \begin{equation}
A=\begin{pmatrix} a_2 & 1 & 0& \dots & & 0\\
a_4 & a_2 & 1 & 0& \dots & 0\\
\vdots & & &\ddots & &\vdots \\
& & & & \ddots & \\
& & & & & 1\\
a_{2m}&  & \dots & a_6 & a_4& a_2
\end{pmatrix}
 \label{symp}
 \end{equation}
 
 Then $\det(\lambda-\Phi)=\det (\lambda^2-A)$ and the  coefficients are again universal polynomials in the $a_i$, providing another basis  for the invariant polynomials.
 \section{Spectral data}\label{spec}
 
 Given a Higgs bundle for a linear group the spectral curve $S$ is defined by the characteristic equation $0=\det(x-\Phi)=x^n+a_1x^{n-1}+\dots+a_n$, where $a_i\in H^0(\Sigma, K^i)$. It is  a divisor of the line bundle $\pi^*K^n$ on the total space of $\pi:K\rightarrow \Sigma$, where  $x$ is the tautological section of $\pi^*K$ on $K$, and since the canonical bundle of a cotangent bundle  is trivial, $K_S\cong \pi^*K^n$ by adjunction. In particular the genus of $S$ is given by $g_S-1=n^2(g-1)$. When $S$ is smooth, the cokernel of $\pi^*\Phi-xI$ on $S$ defines a line bundle $L\pi^*K (=L\otimes \pi^*K)$ and the vector bundle $V$ can be recovered as $V=\pi_*L$, the direct image sheaf. The direct image of $x:L\rightarrow L \pi^*K$ is then the Higgs field $\Phi$. The direct image of the trivial bundle is 
 ${\mathcal O}\oplus K^{-1}\oplus K^{-2}\oplus\dots\oplus K^{-(n-1)}$  \cite{BNR}.   
 
 In general  $\Lambda^nV\cong \Nm(L) K^{-n(n-1)/2}$  \cite{BNR} where $\Nm:\Pic(S)\rightarrow \Pic(\Sigma)$ is the norm map which associates to a divisor $\sum n_ip_i$ on $S$  the divisor $\sum n_i\pi(p_i)$ on $\Sigma$. Thus, to get an $SL(n,\C)$-Higgs bundle we take $L\cong U\pi^*K^{(n-1)/2}$ where $U$ lies in the Prym variety, the kernel of the homomorphism $\Nm:\Pic^0(S)\rightarrow \Pic^0(\Sigma)$. The canonical  section described above is obtained by taking $U$ to be the trivial bundle.
 
 For the symplectic group $Sp(2m,\C)$, the eigenvalues of $\Phi$ occur in pairs $\pm \lambda$ and the equation for the spectral curve has the form $x^{2m}+a_2x^{2m-2}+\dots+a_{2m}=0$. Thus $S$ has an involution $\sigma(x)=-x$. In this case the bundle $U$ must satisfy $\sigma^*U\cong U^*$ \cite{Hit2}. This is the Prym variety for the map to the quotient $S\rightarrow \bar S=S/\sigma$.

 These are the spectral data for the complex groups, next we need to find the restrictions    for the real forms. For the group  $SL(n,\R)$, we need $V$ to be orthogonal and $\Phi$ to be symmetric. This is a fixed point of a holomorphic involution on the Higgs bundle moduli space: $(V,\Phi)\mapsto (V^*,\Phi^T)$. Since the real dimension of $\Hom(\pi_1(\Sigma),SL(n,\R))/SL(n,\R)$ is $2(g-1)\dim SL(n,\R)$, each component of the fixed-point set has complex dimension $\dim {\mathcal M}/2$. The canonical Higgs bundle lies in the fixed-point set as we have seen, and so the $+1$-eigenspace of the action on the tangent space at this point has dimension $\dim {\mathcal M}/2$.
However, $\det (x-\Phi^T)=\det(x-\Phi)$ so the involution acts trivially on the base of the integrable system, which means the action on the tangent space to the fibre is $-1$. Since the fibre is known to be connected, by exponentiation the fixed points correspond to the elements of order $2$ in the Prym variety.  In fact, this  argument holds for any split real form and is dealt with in \cite{LS}.

To see more concretely how the direct image $\pi_*L$ acquires an orthogonal structure when $U^2$ is trivial we use relative duality, with the equivalent condition  $L^2\cong  \pi^*K^{n-1}\cong K_S \pi^*K^*$. 

Relative duality in our situation  states that for any vector bundle $W$ on $S$, 
$(\pi_*W)^*\cong \pi_*( W^*\otimes K_S \pi^*K^*).$ Explicitly, over a regular value $p$ of $\pi$, $$(\pi_*W)_p=\bigoplus_{\pi(u)=p} W_u$$
and at each point $u\in S$ we have the derivative $d\pi_u\in (K_S \pi^*K^*)_u$. Then given $v\in \pi_*(W)_p$, $\xi\in \pi_*(W^*\otimes K_S \pi^*K^*)_p$ the non-degenerate pairing is:
$$\langle v,\xi\rangle= \bigoplus_{\pi(u)=p}\frac{\xi(v)_u}{d\pi_u}.$$
In the neighbourhood of a branch point with the local form $z\mapsto w=z^m$ we  write a local holomorphic section of $\pi_*W$ as $f(z)=b_0(w)+zb_1(w)+\dots +z^{m-1}b_{m-1}(w)$ and then, if  $g(z)=c_0(w)+c_1(w)+\dots +z^{m-1}c_{m-1}(w)$ is a local section of $\pi_*(W^*\otimes K_S \pi^*K^*)$, on approaching the branch point  we have a contribution of 
$$\lim_{z\rightarrow 0}\sum_{i=0}^{m-1} \frac{1}{m\omega^{m-1}z^{m-1}}\langle f(\omega^iz),g(\omega^iz)\rangle=\sum_{k+\ell=m-1} \langle b_k,c_{\ell}\rangle$$
where $\omega$ is a primitive $m$th root of unity. This is still well-defined and non-degenerate.

So, returning to the case $L^2\cong K_S \pi^*K^*$, the duality $V\cong V^*$ is expressed by the quadratic form evaluated on $s\in \pi_*(L)_p$ 
\begin{equation}
(s,s)_p= \bigoplus_{\pi(u)=p}\frac{s^2_u}{d\pi_u}
\label{ortho}
\end{equation}
 which is naturally a sum of squares over regular values. The Higgs field is the direct image of $s\mapsto xs$ and since $(xs)t_u=s(xt)_u$, $\Phi$ is symmetric. 
 
 In the symplectic case $\sigma^*U\cong U^*$ and since  $L\cong U \pi^*K^{(2m-1)/2}$ we have 
 $$\sigma^*L\cong L^* \pi^*K^{2m-1}\cong L^* K_S \pi^*K^*.$$
 Here, we have a non-degenerate bilinear form 
\begin{equation}
\langle v,w\rangle= \bigoplus_{\pi(u)=p}\frac{\sigma^*v(w)}{d\pi_u}
\label{skew}
\end{equation}
 which is skew-symmetric since $d\pi$ has opposite signs at $u,\sigma(u)$. 
 
 We also have the condition that $U^2$ is trivial and so $L^2\cong K_S \pi^*K^*$ . This means that $\sigma^*L\cong L$ and we have an action of $\sigma$ (well-defined modulo $\pm 1$) on $L$. Given an open set $A\subset \Sigma$, $\pi^{-1}(A)\subset S$ is invariant under $\sigma$ and this means that the decomposition of $H^0(\pi^{-1}(A),L)$ into invariant and anti-invariant parts descends to a decomposition $V=W_1\oplus W_2$. Since $\sigma$ interchanges in pairs the $2m$ fibres, $\rk W_1=\rk W_2=m$. 
 
 Now if $s,t\in W_1$ they are represented by local invariant sections of $L$. Then from (\ref{skew}), $(\sigma^*s) t$ is invariant but the denominator $d\pi_u$ is anti-invariant and hence $W_1$ is Lagrangian, so $W_2\cong W_1^*$. We therefore have the required form for $V=W\oplus W^*$. Now since $\sigma(x)=-x$, $\Phi$ interchanges $W$ and $W^*$ and as before $(xs)t=s(xt)$ is symmetric.

 \section{Characteristic classes for $SL(n,\R)$}
 
 In the previous section we saw how the direct image of a line bundle $U$ of order $2$ on the  curve $S$ defines an orthogonal bundle $V$ on $\Sigma$. There are two characteristic classes $w_1(V)$ and $w_2(V)$ but $w_1=0$ if $U$ lies in the Prym variety of $\pi:S\rightarrow \Sigma$. Topologically this means that if we take the  dual homology class $u\in H_1(S,\Z_2)$ of $U\in H^1(S,\Z_2)$ then $\pi_*(u)=0$.

 The second Stiefel-Whitney class is more complicated. 
 
 To discuss the topology of orthogonal bundles on a surface $\Sigma$ we use $KO$-theory, following  \cite{MFA1}. 
 For a compact surface $\Sigma$
$$KO(\Sigma)\cong \Z\oplus H^1(\Sigma,\Z_2)\oplus H^2(\Sigma,\Z_2)$$
where the total Stiefel-Whitney class $w=1+w_1+w_2$ gives an isomorphism of the additive group  $\tilde KO(\Sigma)$ to the multiplicative group $ 1\oplus H^1(\Sigma,\Z_2)\oplus H^2(\Sigma,\Z_2)$.  

Generators are given by holomorphic line bundles $L$ such that $L^2\cong {\mathcal O}$  and the class $\Omega={\mathcal O}_p+{\mathcal O}_p^*-2$ where ${\mathcal O}_p$ is the holomorphic line bundle given by a point $p\in \Sigma$. We write $\alpha(x)\in KO(\Sigma)$ for the class of the line bundle corresponding to $x\in H^1(\Sigma,\Z_2)$. Then $\alpha(0)=1$ and 
$$\alpha(x+y)=\alpha(x)+\alpha(y)-1+(x,y)\Omega$$
where $(x,y)$ is the intersection form. 
This expression is nonlinear as it corresponds to the tensor product of line bundles. Then the isomorphism is determined by the relations $w_1(\alpha(x))=x, w_1(\Omega)=0$ and $ w_2(\Omega)=c_1({\mathcal O}_p)$ mod $2$=$[\Sigma]\in  H^2(\Sigma,\Z_2)$.

For an arbitrary rank $n$ orthogonal bundle $V$ we have 
$$[V]=n-1+\alpha(w_1(V))+w_2(V)\Omega.$$

Now choose a theta characteristic $K^{1/2}$ on $\Sigma$. This is a spin structure and defines a $KO$-orientation. The map to a point gives an invariant which is a spin cobordism characteristic number, an additive homomorphism  $\varphi:KO(\Sigma)\rightarrow \Z_2$. 

Given a holomorphic bundle $V$ with an orthogonal structure,  $\varphi$  is an analytic mod $2$ index $\varphi([V])=\dim H^0(\Sigma,  V\otimes K^{1/2})$ mod $2$.
For $V={\mathcal O}_p+{\mathcal O}^*_p$, then (as in \cite{MFA1}), Riemann-Roch  and Serre duality gives
$$1=\dim H^0(\Sigma, {\mathcal O}_p K^{1/2})-\dim H^0(\Sigma, {\mathcal O}^*_p K^{1/2})$$
and so  $\varphi(V)=1=\varphi(\Omega)$.

\begin{thm} \label{w2} Let $S$ be a smooth curve in the total space of $\pi:K\rightarrow \Sigma$ given by an equation $x^n+a_1x^{n-1}+\dots+a_n=0$ and let  $L$ be a line bundle on $S$ such that  $L^2\cong K_S \pi^*K^*$. Define $V=\pi_*L$,  the direct image bundle given the orthogonal structure described above. 

 Let $K^{1/2}$ be a   theta characteristic on $\Sigma$ with $\varphi_{\Sigma}(1)=0$, and $K_S^{1/2}=L \pi^*K^{1/2}$ the corresponding theta characteristic on $S$. Then
$$w_2(V)=\varphi_S(1)+\varphi_{\Sigma}(\alpha(w_1(V))).$$
\end{thm}

\begin{prf} The class of $V$ in $KO(\Sigma)$ is $n-1+\alpha(w_1(V))+w_2(V)\Omega$ so applying $\varphi_{\Sigma}$,
$$w_2(V)=(n-1)\varphi_{\Sigma}(1)+\varphi_{\Sigma}(\alpha(w_1(V))+\varphi_{\Sigma}(V)$$
and $\varphi_{\Sigma}(1)=0$ by the choice of $K^{1/2}$.

Now, the defining property of the direct image is that $H^0(A,\pi_*L)= H^0(\pi^{-1}(A), L)$ for any open set $A\subset \Sigma$ so with $A=\Sigma$
$$H^0(S,K^{1/2}_S)= H^0(S,L \pi^*K^{1/2})=H^0(\Sigma,\pi_*L\otimes K^{1/2})$$
which gives $\varphi_S(1)=\varphi_{\Sigma}(V)$.
\end{prf}

\begin{rmks}   

\noindent 1. There is an  alternative approach to deriving this formula using the topological definition of the mod 2 index due to Thurston (see \cite{AGH} page 291).  Away from the branch locus $B$, the monodromy of a loop in  $\Sigma\backslash B$  preserves the orthogonal structure on the direct image as a sum of squares and so lies in the group  $B(n)$,  the semi-direct product of the symmetric group $S(n)$ and $(\Z_2)^n$. This group  is a subgroup of $O(n)$ and has a double covering $C(n)\subset Pin(n)$ which is a central extension by $\Z_2$. The authors of \cite{EOP} relate the mod 2 invariant to lifting issues related to  this group. In our context it clearly corresponds to the question of whether the structure group of $V$ lifts to $Spin(n)$, i.e.  whether $w_2(V)=0$ or not. 

\noindent 2. The characteristic class $w_2$ is independent of which spin structure $K^{1/2}$ we choose, (which was why it was convenient to take an even one in the theorem). A better way to formulate this fact is to regard the isomorphism $L^2\cong K_S \pi^*K^*$ as a $KO$-theory orientation on the {\it map} $\pi:S\rightarrow \Sigma$. There is then a push-forward map $\pi_{!}:KO(S)\rightarrow KO(\Sigma)$ and $[V]= \pi_{!}(1)$.
\end{rmks}
 
 Given the formula in Theorem \ref{w2} we need to determine how many points of order $2$ give $w_2(V)=0$. The interpretation of $w_2$ via $\varphi_S$ tells us that this is an affine  quadratic function $\psi$ on $\PP[2]$, the elements of order $2$ in the Prym variety.  
 Choosing an origin such that $\psi(0)=0$ this means that $\psi(x+y)=\psi(x)+\psi(y)+(x,y)$ for a bilinear form $(x,y)$. Here it is  the intersection form on $H^1(S,\Z_2)$ restricted to $\PP[2]$.  
 
 The quadratic functions $xy$ and $x^2+xy+y^2$ on $(\Z_2)^2$ have the same bilinear form but are not equivalent and for quadratic functions associated to a nondegenerate bilinear form there are two canonical forms: a sum of $k$ terms of the form $xy$ or a sum of $(k-1)$ such terms plus a single term $x^2+xy +y^2$. They are distinguished by their  Arf invariant $\in \Z_2$ which is zero in the first case and $1$ in the second.  When the invariant is $0$, $\psi$ has $2^{k-1}(2^k+1)$ zeros and otherwise $2^{k-1}(2^k-1)$. This interpretation shows that  the invariant is  independent of any choice of origin. It is additive under orthogonal direct sum  since  $x^2+xy+y^2+u^2+uv+v^2=(x+u)(x+y+u)+(y+v)(y+u+v).$

 Note in what follows that $\pi^*:H^1(\Sigma,\Z_2)\rightarrow H^1(S,\Z_2)$ is injective: indeed, (as in \cite{BNR}), given a degree zero line bundle $L$ on $\Sigma$ with $\pi^*L$ trivial, we have a section of $\pi_* \pi^*L=L\oplus LK^{-1}\oplus\dots$. Since $L$ has degree $0$ and  $g>1$ this must be  a section of $L$  which is therefore  trivial. 
 
 \begin{prp} In the context of Theorem \ref{w2}, when $n$ is odd  there are $2^{2p-1}+2^{p-1}$ choices of $L$  which give $w_1(V) = w_2(V)=0$; and when $n=2m$ there are $2^{2p-1}+(-1)^{m(g-1)}2^{p+g-1}$ choices, where  $p=(g-1)(n^2-1)$.  \end{prp}
  
  \begin{prf} We have already observed that  $w_1(V)=0$ if $U$ lies in the  Prym variety $\PP$.   
  
  \noindent 1. First consider the case where  $n=2m+1$ is odd. 
  
  The Prym variety has polarization $(1,1,1,\dots, n,n,\dots n)$ with $g$ copies of $n$ (see \cite{BNR})   
so since $n$ is odd the intersection matrix mod $2$ is  non-degenerate. Moreover, since $\Nm \pi^*(x)=(2m+1)x=x$ if $2x=0$ then $H^1(S,\Z_2)$ is an  orthogonal direct sum  $\pi^*H^1(\Sigma,\Z_2)\oplus \PP[2]$.  
  
  From Theorem \ref{w2}  we have  $w_2(V)=\varphi_S(1)$ as a function of theta characteristics of the form $K_S^{1/2}=L K^{1/2}$. By \cite{MFA1}  for all choices of $K_S^{1/2}$ the Arf invariant is $0$.
   But  the   invariant  is additive under orthogonal direct sum, and if we take 
    $K^{1/2}_S=\pi^*(U K^{m+1/2})$ for $U\in H^1(\Sigma,\Z_2)$, then taking the direct image,
\begin{eqnarray*}
\dim H^0(S,K^{1/2}_S)&=&\dim H^0(\Sigma,U\otimes (K^{-m+1/2}\oplus \dots\oplus K^{m+1/2}))\\
&=&\dim H^0(\Sigma, U\otimes K^{1/2})+2(g-1)+\dots +2m(g-1)
\end{eqnarray*}
which is $\dim H^0(\Sigma, U\otimes K^{1/2})$ mod $2$.   This is the standard quadratic function for the Riemann surface $\Sigma$ so the Arf invariant  is zero and by additivity so is the invariant on $\PP[2]$. 

We have $\dim H^0(\Sigma, K^{1/2})=0$ mod $2$ by choice, the origin,   and so  it follows that  there are $2^{p-1}(2^p+1)$ zeros where $p=(g-1)\dim SL(n)=(g-1)(n^2-1)$ is the dimension of the Prym variety. 

\noindent 2. Now assume $n=2m$.

 In this case $\Nm \pi^*(x)=2mx=0$ and $\pi^*H^1(\Sigma,\Z_2)$ lies inside  $\PP[2]$, as the degeneracy subspace of the intersection form.   Then $K^{1/2}_S=\pi^*(U K^{m})$ for $U\in H^1(\Sigma,\Z_2)$ we have 
   \begin{eqnarray*}
\dim H^0(S,K^{1/2}_S)&=&\dim H^0(\Sigma,U\otimes (K^{-m+1}\oplus \dots\oplus K^{m}))\\
&=&\dim H^0(\Sigma, U)+H^0(\Sigma,U K)+3(g-1)+\dots +(2m-1)(g-1)
\end{eqnarray*}
  which by Riemann-Roch and Serre duality is $m(g-1) $ mod $2$. 
  
  From the proof of  Theorem \ref{w2} in this case $w_2(V)=\varphi_S(1)+2m\varphi_{\Sigma}(1)=\varphi_S(1)$ independently of the choice of $K^{1/2}$. In particular this means that $\varphi_S\vert_{\PP[2]}$ is invariant under the action  of $U\in \pi^*H^1(\Sigma,\Z_2)$ and hence 
 for $y\in \PP[2]$ and $x\in \pi^*H^1(\Sigma,\Z_2)$, $\psi(x+y)=\psi(y)$, and so $\psi(x)=0$ and $(x,y)=0$. 

 Choose a transverse $2p-2g$-dimensional subspace $X$ to $\pi^*H^1(\Sigma,\Z_2)$ and consider the quadratic function $\psi$ restricted to $X$. Then from the canonical form there is a basis  $y_i,z_i$ of  $X$ such that the function is $\sum_{i=1}^{p-g} a_ib_i$ or $\sum_{i=2}^{p-g} a_ib_i+a_1^2+a_1b_1+b_1^2.$

Take  a basis $x_1,\dots x_{2g}$ for $\pi^*H^1(\Sigma,\Z_2)$, then by non-degeneracy of the intersection form on $H^1(S,\Z_2)$ there are elements $w_i$ such that $(x_i,w_j)=\delta_{ij}$.

Define $$\tilde w_1=w_1+\sum_i(w_1,y_i)z_i+\sum_i(w_1,z_i)y_i$$
then  $\tilde w_1$ is orthogonal to $X$. Since each $x_i$ is orthogonal to $\PP[2]$ the $2$-dimensional space spanned by $x_1,\tilde w_1$ is orthogonal to $X$. So since $\psi(x_1)=0$ 
$$\psi(a\tilde w_1+bx_1)=a^2\psi(\tilde w_1)+b^2\psi(x_1)+ab=a(a\psi(\tilde w_1)+b)$$
and hence on the space spanned by $X$ and these two vectors we are adding an $xy$ term, which means we have the same  Arf invariant. By induction so does the full space. 

Now, as we showed above, for $U\in \pi^*H^1(\Sigma,\Z_2)$ we have $\dim H^0(S,K^{1/2}_S)=m(g-1)$ mod $2$, so the number of zeros on $X$ is 
$2^{p-g-1}(2^{p-g}+(-1)^{m(g-1)})$. Acting by $U\in \pi^*H^1(\Sigma,\Z_2)$ gives all of $\PP[2]$ and hence 
the total number of zeros is 
\begin{equation}
2^{2g}\times 2^{p-g-1}(2^{p-g}+(-1)^{m(g-1)})=2^{2p-1}+(-1)^{m(g-1)}2^{p+g-1}
\label{numbers}
\end{equation}
\end{prf}

  \section{Characteristic classes for $Sp(2m,\R)$}\label{sympl}
  
  We have already seen in Section \ref{spec} that a Higgs bundle for the group $Sp(2m,\R)$ is obtained from the spectral curve $S$ defined by $x^{2m}+a_2x^{2m-2}+\dots +a_{2m}=0$ as the direct image $V=W\oplus W^*$ of a line bundle $L\cong U \pi^*K^{(2m-1)/2}$ such that $\sigma^*U\cong U$ and $U^2$ is trivial. Moreover $W$ is the direct image of the invariant sections of $L$. The characteristic class here is $c_1(W)$, and we need to evaluate this as a function  on $\PP[2]$, the points of order $2$ on the Prym variety of  $p:S\rightarrow S/\sigma=\bar S$. 
  
  Note first that, as in the case of $SL(2m,\R)$, $p^*H^1(\bar S,\Z_2)$ lies in the Prym variety. The map $p^*$  is also injective. This is a similar argument to the one above. In this case  the invariant and anti-invariant parts decompose the direct image of $p^*M$ as $M\oplus M'$ where $M'$ has negative degree. So if $p^*M$ has a section, so does $M$.   The dimension of the $\Z_2$-vector space $\PP[2]/p^*H^1(\bar S,\Z_2)$ is therefore $2(g_S-2g_{\bar S})$, but by Riemann-Hurwitz, since $S\rightarrow \bar S $ has $4m(g-1)$ branch points, 
  $2-2g_S=2(2-2g_{\bar S})-4m(g-1)$ and so 
  $$\dim P[2]/p^*H^1(\bar S,\Z_2)=4m(g-1)-2.$$

  We now use the condition $\sigma^*U\cong U$, so that the involution lifts to the line bundle $U$. There are two lifts $\pm \sigma$ but fix attention on one for the moment. Following \cite{LS1},  we consider  the action $\pm 1$ of $\sigma$ on the fibre of $U$ at a fixed point.  
   
 \begin{prp} Suppose the action is $-1$ at $\ell$ fixed points, then $c_1(W)=-\ell/2+m(g-1)$.
  \end{prp}
 
 \begin{prf} The fixed point set of $\sigma(x)=-x$ is the intersection of the zero section of $K$ with $S$. Setting $x=0$ in the equation $x^{2m}+a_2x^{2m-2}+\dots +a_{2m}=0$, these points are the images of the $4m(g-1)$ zeros of $a_m\in H^0(\Sigma,K^{2m})$ under the zero section. The action is $-1$ at $\ell$ of these points. 
 
  Choose  a line bundle  $M$ on $\Sigma$ of large enough degree that  the higher cohomology groups vanish and then 
 applying the holomorphic Lefschetz formula  \cite{AB} we obtain
 $$\dim H^0(S,L \pi^*M)^+-\dim H^0(S,L \pi^*M)^-=\frac{1}{2}(-\ell+(4m(g-1)-\ell))$$
 where the superscript denotes the $\pm 1$ eigenspace under the action of $\sigma$.  Riemann-Roch gives 
  $$\dim H^0(S,L \pi^*M)^++\dim H^0(S,L \pi^*M)^-=\dim H^0(\Sigma, V\otimes M)=2m(1-g+c_1(M))$$ since $V$ is symplectic and $\deg V=0$.  Hence
  $$\dim H^0(S,L\otimes \pi^*M)^+=-\frac{\ell}{2}+mc_1(M)$$

  But by the definition of $W$  this is 
 $\dim H^0(\Sigma,W\otimes M)$ and Riemann-Roch and the vanishing of $H^1$ gives the value $m(1-g)+c_1(W)+mc_1(M)$ and so
 \begin{equation}
 c_1(W)=-\frac{\ell}{2}+m(g-1).
 \label{degree}
 \end{equation}

  \end{prf} 
  \begin{rmks}
  
   \noindent 1. Since $0\le \ell\le 4m(g-1)$ we have $\vert c_1(W)\vert \le m(g-1)$ which is the Milnor-Wood inequality for this group \cite{GGM}. 
  
  \noindent 2. Taking the action $-\sigma$ instead of $\sigma$ changes $\ell$ in the formula to $4m(g-1)-\ell$ and $c_1(W)$ to $-c_1(W)$, and then the roles of $W$ and $W^*$ are interchanged.
  \end{rmks}
  
  The formula (\ref{degree}) above clearly requires $\ell$ to be even, but there is a  reason for this.  If $\ell=0$ the action is trivial at all fixed points, so  $U$ is the pull-back of a flat line bundle of order $2$ on the quotient $\bar S$.  In  general, let  $B$ denote  the subset of $\ell$ points in the branch locus of $S\rightarrow \bar S$, then the line bundle  corresponds to a flat line bundle on $\bar S\backslash B$  where the local holonomy around each $b\in B$ is $-1$. The global holonomy defines a homomorphism  $\rho:\pi_1(\Sigma\backslash B)\rightarrow \Z_2$ 
where the  fundamental group has generators $A_i,B_i$, $1\le i\le g$ and $\delta_j$, $1\le j\le N$, each $\delta_j$  defining  a loop around a branch point. These  satisfy the relation
$$\prod_i[A_i,B_i]\prod_j\delta_j=1$$
but then, with values in the abelian group $\Z_2$, we must have $\prod_j\delta_j=1$ and hence an even number of $-1$ terms.

 This interpretation helps to understand which of the $2^{2p}$ elements in the Prym variety yield a given characteristic class. Let $Z$ be the $4m(g-1)$-element set of zeros of the section $a_{2m}$ of $K^{2m}$ and  let $C(Z)$ be the space of $\Z_2$-valued functions on $Z$, and $C_0(Z)$ the subspace of those whose integral is zero, i.e. take the value $1$ an even number of times. The constant function $1$ lies in  $C_0(Z)$ and let $H(Z)$ be the $(4m(g-1)-2)$-dimensional quotient. 
 
 For a line bundle $U\in \PP[2]$ let $A$ be the subset of $Z$ on which the action of $\sigma$ is $-1$. As noted above, $A$ has an even number of elements  and so its characteristic  function $\chi_A$ lies in $ C_0(Z)$. Define $f(U)\in H(Z)$ to be its  
 equivalence class.  Since we take the quotient by the constant function $1$, this is independent of the choice of lift of $\sigma$. 
 
 \begin{prp} \label{hz} The homomorphism $f$ from $\PP[2]$ to $H(Z)$ is surjective and has kernel $p^*H^1(\bar S,\Z_2)\subset \PP[2]$.
  \end{prp} 
  
  \begin{prf}
 Given {\it any} subset of the zeros of $a_{2m}$  with an even number of elements    we can choose $\delta_j=-1$ as above and get a flat line bundle $U\in H^1(S,\Z_2)$ and  an action of $\sigma$ which acts as $-1$ at those points, so the homomorphism is surjective.  
 The kernel consists of line bundles with trivial action at all fixed points and these are precisely those pulled back from the quotient $\bar S$.
  \end{prf}
  Since $H(Z)$ and $\PP[2]/p^*H^1(\bar S,\Z_2)$ have the same number of elements we see from the Proposition that they are isomorphic. 
  
We can now count the points in $\PP[2]$ with fixed characteristic class. This is, from Equation (\ref{degree}),  $c_1(W)=m(g-1)-\ell/2$ and  involves a choice between $W$ and $W^*$  or equivalently a choice of lifting of the involution $\sigma$, or the subset $A\subset Z$ with an even number $\ell$ of elements.  Thus 
the number of such elements in $\PP[2]$ is, from Proposition \ref{hz}, 
$${4m(g-1)\choose \ell}\times 2^{2q}$$
where $q=g_{\bar S}=(2m^2-m)(g-1)+1$.
 
\begin{rmks} 

\noindent 1. The $\Z_2$-vector space $H(Z)$ is a representation of the symmetric group $S(4m(g-1))$, the permutations of the $4m(g-1)$ branch points, and  we may describe the above result by saying that the characteristic class is determined by the orbit of the symmetric group on this space. In the case of $n=2$ this picture was derived in \cite{LS0} via the monodromy action of the family of abelian varieties.

 \noindent 2. For $n=2$ the two groups coincide, so we may use the formula above to compare with the $SL(2,\R)$ case. Here $S$ is a double cover of $\Sigma$ and so $\bar S=\Sigma$, hence $H(Z)=\PP[2]/\pi^*H^1(\Sigma,\Z_2)$. Equation (\ref{degree})  
gives  $c_1(W)=-{\ell}/{2}+(g-1)$ and $w_2=c_1(W)$ mod $2$, but here we don't distinguish between $W$ and $W^*$ so the number with $w_2=0$ is 
$$\frac{1}{2}\sum_{\ell\equiv (2g-2) {\mathrm{mod}}  \, 4}{4(g-1)\choose \ell}\times 2^{2g}$$
If $(g-1)$ is even this is 
$$2^{2g-3}((1+1)^{4(g-1)}+(1+i)^{4(g-1)}+(1-1)^{4(g-1)}+(1-i)^{4(g-1)})$$
and  if $(g-1)$ is odd
$$2^{2g-3}((1+1)^{4(g-1)}+(1+i)^{4(g-1)}-(1-1)^{4(g-1)}-(1-i)^{4(g-1)}).$$
Using $e^{i\pi/4}=(1+i)/
\sqrt{2}$ this gives $2^{6g-7}+2^{4g-4}$ or $2^{6g-7}-2^{4g-4}$ which checks with (\ref{numbers}). 
\end{rmks}

\section{Mirror symmetry}
A hyperk\"ahler manifold has three complex structures $I,J,K$ but also three symplectic forms, the corresponding K\"ahler forms $\omega_1,\omega_2,\omega_3$. A {\it brane} for a complex manifold (a B-brane) is roughly speaking a holomorphic bundle over a complex submanifold and for a symplectic manifold (an A-brane) it is a flat vector bundle over a Lagrangian submanifold. For a hyperk\"ahler manifold a BAA-brane is a B-brane for the complex structure $I$ and an A-brane for the symplectic structures $\omega_2,\omega_3$. The trivial bundle over a component of the $G^r$-character variety is an example. As the fixed point set of an antiholomorphic isometry for complex structures $J,K$ it is Lagrangian for the K\"ahler forms  $\omega_2,\omega_3$, and it is holomorphic with respect to $I$. 

Mirror symmetry should transform a BAA-brane on a hyperk\"ahler manifold $M$  to a BBB-brane on its mirror $\hat M$ \cite{KW} (I am indebted to Sergei Gukov for this information). A BBB-brane is a holomorphic bundle over a complex submanifold  with respect to all complex structures $I,J$ and $K$, equivalently  a hyperk\"ahler submanifold with a hyperholomorphic bundle over it. A hyperholomorphic bundle is a bundle with connection whose curvature is of type $(1,1)$ with respect to all complex structures. Such connections (generalizations of instantons in four dimensions) are quite rare and so it is intriguing to seek such an object as the mirror of a $G^r$-character variety. We shall attempt this now for the group $G^r=U(m,m)\subset GL(2m,\C)$.

The Higgs bundle and spectral data description of $U(m,m)$ will be  familiar from our previous discussion of $Sp(2m,\R)$. Details can be found in \cite{LS1}. The Higgs bundle is of the form $V=W_1\oplus W_2$  and the Higgs field $\Phi$ is off-diagonal: $\Phi(w_1,w_2)=   (\beta(w_2),\gamma(w_1))$. There is a characteristic class $c_1(W_1)$, and to keep the link with flat connections we need $c_1(V)=0$ and so $c_1(W_1)=-c_1(W_2)\in \Z$. The spectral curve $S$ has the form $x^{2m}+a_2x^{2m-2}+\dots+a_{2m}=0$ and hence an involution $\sigma(x)=-x$ and the spectral data consist of taking  a line bundle $L$ on $S$ such that $\sigma^*L\cong L$. As in Section \ref{sympl}, the characteristic class is determined by the number of points on $x=0$ at which the lifted action of $\sigma$ is $-1$. The difference here with the $Sp(2m,\R)$-case is that the fibre is not discrete but is instead the disjoint union of a finite number of abelian varieties. In fact if $L_1,L_2$ are two line bundles with the same subset of fixed points at which $\sigma$ acts as $-1$, then the action on $L^*_1L_2$ is trivial and so it is pulled back from the quotient $\bar S$. Thus the fibre is isomorphic to  the disjoint union of  $N$ copies of  $\Pic^0(\bar S)$ where $N=2^{4m(g-1)-1}$ is the number of  subsets of the zero set $Z$ of $a_{2m}$ with an even number of elements. 

 For a Calabi-Yau manifold with a special Lagrangian fibration mirror symmetry is effected via the Strominger-Yau-Zaslow approach of replacing each nonsingular torus fibre by its dual, and hoping that it can be extended over the discriminant locus in the base. The Higgs bundle integrable system fits into this framework as first investigated in \cite{HT}. As in \cite{Hit4} for certain cases and \cite{DP} in general, it corresponds to replacing the group $G^c$ by its Langlands dual group $^LG^c$. 

 We now consider the structure on the dual fibration relevant for $U(m,m)\subset GL(2m,\C)$. The abelian variety for $GL(2m,\C)$ is the Jacobian, or  $\Pic^0(S)$, and  since Jacobians are self-dual, and the norm map is the adjoint of the pull-back,  dualizing the inclusion $\Pic^0(\bar S)\subset \Pic^0(S)$ gives 
$$0\rightarrow \PP(S,\bar S)\rightarrow \Pic^0(S)\rightarrow  \Pic^0(\bar S)\rightarrow 0$$
and the Prym variety $\PP(S,\bar S)$ is a distinguished subvariety   of $\Pic^0(S)$. In terms of duality it parametrizes line bundles on $\Pic^0(S)$ which are trivial on $\Pic^0(\bar S)$.

If mirror symmetry is to work as predicted this family of abelian varieties over the space of polynomials  $x^{2m}+a_2x^{2m-2}+\dots+a_{2m}$ should extend to a hyperk\"ahler submanifold of the Higgs bundle moduli space for the Langlands dual of $GL(2m,\C)$, which is again $GL(2m,\C)$. But in Section \ref{spec} we saw that the $Sp(2m,\C)$-moduli space  had an integrable system over this base  whose generic fibre was $\PP(S,\bar S)$. The inclusion of a group gives a hyperk\"ahler subspace of the moduli space so the symplectic Higgs bundles form a hyperk\"ahler subspace of the moduli space of $GL(2m,\C)$-Higgs bundles. We therefore have the first requirement of the mirror -- the hyperk\"ahler support for a hyperholomorphic bundle.

The remaining task is to find a hyperholomorphic vector bundle over the $Sp(2m,\C)$-moduli space, or rather several, one for each characteristic class. There are relatively few constructions of such bundles but there is one which involves a Dirac-type operator, and which we describe next.  More information may be found in \cite{Hit5},\cite{Bon}. 

For each $(A,\Phi)$ satisfying the Higgs bundle equations for a compact group $G$, we take a vector bundle $V$ associated to the principal $G$-bundle via a representation of $G$ and define an elliptic operator $D^*:V\otimes (K\oplus \bar K)\rightarrow V\otimes (K\bar K\oplus K\bar K)$ by 
$$D^*=\begin{pmatrix} \bar\partial_A & \Phi\\
\Phi^*& \partial_A
\end{pmatrix}.$$
 The equation $F_A+[\Phi,\Phi^*]=0$ yields a vanishing theorem for irreducible connections and the index theorem gives $\dim \ker D^*=(2g-2)\rk V$. For the adjoint representation the null space can be viewed as  the tangent space of the moduli space. 
 
 Given a universal bundle over ${\mathcal M}\times \Sigma$, the family of null-spaces  defines a rank  $(2g-2)\rk V$ vector bundle on ${\mathcal M}$, and since  $D^*$ acts on one-forms with values in a Hermitian bundle $V$, there is a conformally-invariant ${\mathcal L}^2$ inner product  which defines by projection a connection on this bundle. It turns out that this connection is hyperholomorphic: for the adjoint representation  it is the Levi-Civita connection. We shall leave till later the issue of the existence of a universal bundle -- locally these exist and connections are locally determined. 
 
 The null-space of $D^*$ can be viewed in different ways according to the complex structures $I,J,K$. For $J$ the operator $D^*$ is the Hodge operator for the de Rham complex of the flat connection $\nabla_A+\Phi +\Phi^*$; for $I$ it is the Hodge operator for the total differential $\bar\partial \pm \Phi$ in the double complex 
 $$\Omega^{0,\ast}(V)\stackrel{\Phi}\rightarrow \Omega^{0,\ast}(V).$$
 In this latter case $\ker D^*$ is identified with the hypercohomology group $\K^1$. Trading the Dolbeault viewpoint for the \v{C}ech approach, it is the hypercohomology for the complex of sheaves
 $${\mathcal O}(V)\stackrel{\Phi}\rightarrow {\mathcal O}(V\otimes K).$$   
 
 Now take $V$ to be associated to the defining representation of $Sp(m)$, and $(A,\Phi)$ to lie in a generic fibre of the integrable system. Then $\Phi:V\rightarrow V\otimes K$ is generically an isomorphism and the spectral sequence for the complex of sheaves identifies the hypercohomology group $\K^1$ with sections of a sheaf supported on the zero set of $\det \Phi$. For a smooth spectral curve, $\det \Phi=a_{2m}$ has $4m(g-1)$ distinct zeros $Z$ and we have 
 $$\K^1\cong \bigoplus_{z\in Z}\coker \Phi_z.$$
From the spectral data, $L \pi^*K=\coker (x-\Phi)$, so identifying $Z$ with the intersection of the spectral curve $S$ with the zero section $x=0$ of $K$ we have 
$$\K^1\cong \bigoplus_{z\in Z}(L\pi^* K)_z$$
(note  in particular the dimension checks with the index theory calculation). 

To summarize, what we have here is a hyperholomorphic bundle ${\mathbf V}$ over the $Sp(2m,\C)$-moduli space whose fibre at a point defined by a nonsingular spectral curve $S$ and line bundle $L$ is given by $ \oplus_{z\in Z}(L \pi^*K)_z$. 

Fix the spectral curve $S$, then  the components of the $U(m,m)$-character variety with a fixed characteristic class correspond to the subsets of $\ell=2k$ elements in $Z$. With the same spectral curve,   the corresponding  fibre of the  $Sp(2m,\C)$-moduli space  consists of  line bundles $L$ in the Prym variety and 
$$\Lambda^{2k}{\mathbf V}\cong \bigoplus_{\{z_1,\dots, z_{2k}\}\subset  Z}(L\pi^*K)_{z_1}(L\pi^*K)_{z_2}\dots(L\pi^*K)_{z_{2k}}$$
is a sum over all such subsets of $2k$ elements. This bundle, with its induced hyperholomorphic connection, therefore seems a natural choice for the mirror: as a direct sum over components of the fibre it is analogous to the Fourier-Mukai transform yet it is well-defined on the whole moduli space apart from the issue of the universal bundle, which we consider next. 

The $Sp(2m,\C)$-moduli space has no universal bundle:  the obstruction lies in $H^2({\mathcal M},\Z_2)$. To be more concrete, for an open  covering $\{U_{\alpha}\}$ of ${\mathcal M}$ there is a local universal bundle ${\mathcal V}_{\alpha}$  on $U_{\alpha}\times \Sigma$ and on $U_{\alpha}\cap U_{\beta}$ there is a line bundle ${\mathcal L}_{\alpha\beta}$ of order $2$  such that ${\mathcal V}_{\beta}\cong {\mathcal V}_{\alpha}\otimes {\mathcal L}_{\alpha\beta}$ with compatibility conditions on the isomorphisms. This describes the gerbe which is the obstruction. The $D^*$ null-spaces define bundles ${\mathbf V}_{\alpha}$ over $U_{\alpha}$ with a hyperholomorphic connection and these are related on the intersection  $U_{\alpha}\cap U_{\beta}$ by the flat line bundle ${\mathcal L}_{\alpha\beta}$. However, the even exterior power $\Lambda^{2k}{\mathbf V}_{\alpha}$ is insensitive to this ambiguity and so is well-defined globally.

\begin{rmks}

\noindent 1. The operator $D^*$ may be regarded as a quaternionic operator (see \cite{Hit5}) with coefficients in the bundle $V$ whose structure group is $Sp(m)$ hence also quaternionic. This means the null-space has a real structure and a Hermitian structure, equivalently an orthogonal structure. In particular  we have an isomorphism as bundles with connection $\Lambda^{2k}{\mathbf V}\cong \Lambda^{4m(g-1)-2k}{\mathbf V}$. As we have seen it is only the choice of lifting of the involution $\sigma$ that distinguishes a $2k$-element subset of $Z$ and its complement, so this is expected.

\noindent 2. From the point of view of the spectral data there is a natural orthogonal structure which is almost certainly the same as the above differential-geometric description. At each point  $z\in Z$
we have $a_{2m}(z)=0$ and the derivative at a simple zero defines $da_{2m}(z)\in K^{2m+1}_z$. But $L\cong U \pi^*K^{(2m-1)/2}$ and $\sigma^*U\cong U^*$, so at the fixed point $z$ of $\sigma$ we have a non-zero vector $u_z$ in $U^{-2}_z$. Given $s\in (L\pi^*K)_z$ we can define 
$s^2u_z/da_{2m}(z)\in \C$ and summing these get a non-degenerate quadratic form on ${\mathbf V}$.

\noindent 3. The cases $k=0$ or $k=2m(g-1)$ correspond to the maximal absolute value of the characteristic class allowed and here the hyperholomorphic bundle is the trivial line bundle $\Lambda^0{\mathbf V}_{\alpha}$. Maximal representations play a special role in the study of character varieties (see e.g. \cite{Burg}). 
\end{rmks}


\begin{thebibliography}{11}
   \bibitem{AGH}
   E.Arbarello, M.Cornalba, P.A. Griffiths \& J.Harris, ``Geometry of algebraic curves".  Vol. I.  Springer-Verlag, New York, (1985).    
   %
   \bibitem{AB}
   M.F.Atiyah \& R.Bott,
   A Lefschetz fixed point formula for elliptic complexes II. Applications, {\it Ann. of Math.} {\bf 88} (1968), 451--491. 
   %
   \bibitem{MFA1}
M.F.Atiyah, Riemann surfaces and spin structures, {\it  Ann. Sci. \'Ecole Norm. Sup.} {\bf 4} (1971), 47Ð- 62. 
%
\bibitem{BNR}
A.Beauville, M.S.Narasimhan \& S.Ramanan, 
 Spectral curves and the generalised theta divisor,  
{\it J. Reine Angew. Math.}  {\bf 398}  (1989), 169--179.
%
\bibitem{Bon}
 J.Bonsdorff, 
 Autodual connection in the Fourier transform of a Higgs bundle, {\it  Asian J. Math.} {\bf  14} (2010,) 153Ð-173.
 %
 \bibitem{Burg}
 M.Burger, A.Iozzi \& A.Wienhard,  Surface group representations with maximal Toledo invariant, {\it  Ann. of Math.} {\bf 172} (2010), 517Ð-566.
 %
\bibitem{DP}
R.Donagi \& T.Pantev,  Langlands duality for Hitchin systems, {\it  Invent. Math.} {\bf  189} (2012) 653Ð-735.
%
\bibitem{EOP}
A.Eskin, A.Okounkov \& R.Pandharipande, 
The theta characteristic of a branched covering, 
{\it Advances in Math.}  {\bf 217}  (2008), 873--888.
%
%
\bibitem{GGM}
O.Garcia-Prada, P.Gothen \& I. Mundet i Riera, Higgs bundles and surface group representations in the real symplectic group, {\it J.Topol.} {\bf 6} (2013), 64--118. 
%
\bibitem{HT}
T.Hausel \& M.Thaddeus, Mirror symmetry, Langlands duality and Hitchin systems, {\it Inventiones Mathematicae}, {\bf 153} (2003), 197--229.
%
  \bibitem{Hirz1}
    F.Hirzebruch,
Problems on differentiable and complex manifolds,  
{\it Ann. of Math.} {\bf  60} (1954), 213--236.

%
\bibitem{Hit1}
N.J.Hitchin,
The self-duality equations on a Riemann surface, {\it Proc.London Math.Soc.} {\bf 55}
(1987), 59Ð-126.
%
\bibitem{Hit2}
N.J.Hitchin,
Stable bundles and integrable systems, {\it Duke Math.J.} {\bf 54} (1987),
 91--114.
 %
 \bibitem{Hit3}
 N.J.Hitchin,
 Lie groups and Teichm\"uller space, {\it Topology} {\bf 31} (1992), 
449 -- 473.
%
\bibitem{Hit5}
 N.J.Hitchin,
{\it The Dirac operator}, in ``Invitations to Geometry and Topology'', M.Bridson and S.Salamon (eds.),  Oxford Graduate Texts in Mathematics, Oxford University Press, Oxford (2002), 208 -- 232.
%
 \bibitem{Hit4}
 N.J.Hitchin,
Langlands duality and $G_2$ spectral curves, {\it Quart. J. Math. Oxford}. {\bf 58} (2007), 319--344.
%
\bibitem{KW}
A.Kapustin \& E.Witten, 
Electric-magnetic duality and the geometric Langlands program, 
{\it Commun. Number Theory Phys.} {\bf  1} (2007), 1--236. 
%
\bibitem{LS0}
L.P.Schaposnik,
Monodromy of the $SL_2$ Hitchin fibration, arXiv 111.2550, {\it International Journal of Mathematics} {\bf 24} (2013).
%
 \bibitem{LS}
L.P.Schaposnik,
Spectral data for $G$-Higgs bundles, D.Phil thesis (Oxford), 2013. 
%
 \bibitem{LS1}
L.P.Schaposnik,
Spectral data for $U(m,m)$-Higgs bundles, arXiv 1307.4419.
%
 \bibitem{Sim}
  C.Simpson, 
   Higgs bundles and local systems, {\it Inst.Hautes \'Etudes Sci.Publ.Math.} {\bf 75} (1992), 5--95.
  %
  \bibitem{Tr} 
   W.Trench,  On the eigenvalue problem for Toeplitz band matrices, {\it Linear Algebra Appl.} {\bf  64} (1985), 199--214.
  \end{thebibliography}
 \end{document}